\documentclass[12pt,reqno]{amsart}
\usepackage{amsmath , amssymb}
\usepackage[mathscr]{eucal}

\newtheorem{thm}{Theorem}

\newtheorem{lemma}[thm]{Lemma}

\newtheorem{obs}{Observation}

\newcommand{\cur}{\mathcal}

\def\star(#1,#2){{\bf (\textasteriskcentered)$_{#1,#2}$}}

\def\Fset{\mathcal{F}}

\def\P{\mathcal{P}}
\def\Q{\mathcal{Q}}
\def\T{\mathcal{T}}

\def\R{\mathcal{R}}

\def\G{\mathcal{G}}

\def\K{\mathcal{K}}
\def\L{\mathcal{L}}
\def\M{\mathcal{M}}

\def\N{\mathbb{N}}

\def\C{\mathcal{C}}

\def\le{\leqslant}
\def\ge{\geqslant}

\begin{document}

\title[Hereditary properties: posets and oriented graphs]{Hereditary properties of combinatorial structures:\\ posets and oriented graphs}

\author{J\'ozsef Balogh}
\address{Department of Mathematics\\ University of Illinois\\ 1409 W. Green Street\\ Urbana, IL 61801} \email{jobal@math.uiuc.edu}

\author{B\'ela Bollob\'as}
\address{Trinity College\\ Cambridge CB2 1TQ\\ England \\ and\\ Department of Mathematical Sciences\\ The University of Memphis\\ Memphis, TN 38152} \email{B.Bollobas@dpmms.cam.ac.uk}

\author{Robert Morris}
\address{Instituto Nacional de Matem\'atica Pura e Aplicada, Estrada Dona Castorina, 110, Jardim Bot\^{a}nico, Rio de Janeiro, Brazil (Work done while at The University of Memphis.)} \email{rdmorrs1@impa.br}\thanks{The first author was supported during this research by OTKA grant T049398 and NSF grants DMS-0302804,  DMS-0603769 and DMS 0600303, UIUC Campus Research Board 06139, and OTKA 049398, the second by ITR grant CCR-0225610 and ARO grant W911NF-06-1-0076, and the third by a Van Vleet Memorial Doctoral Fellowship.}


\begin{abstract}
A hereditary property of combinatorial structures is a collection of structures (e.g. graphs, posets) which is closed under isomorphism, closed under taking induced substructures (e.g. induced subgraphs), and contains arbitrarily large structures. Given a property $\P$, we write $\P_n$ for the collection of distinct (i.e., non-isomorphic) structures in a property $\P$ with $n$ vertices, and call the function $n \mapsto |\P_n|$ the \textit{speed} (or unlabelled speed) of $\P$. Also, we write $\P^n$ for the collection of distinct \emph{labelled} structures in $\P$ with vertices labelled $1,\ldots,n$, and call the function $n \mapsto |\P^n|$ the \textit{labelled speed} of $\P$.

The possible labelled speeds of a hereditary property of graphs have been extensively studied, and the aim of this paper is to investigate the possible speeds of other combinatorial structures, namely posets and oriented graphs. More precisely, we show that (for sufficiently large $n$), the labelled speed of a hereditary property of posets is either $1$, or exactly a polynomial, or at least $2^n - 1$. We also show that there is an initial jump in the possible unlabelled speeds of hereditary properties of posets, tournaments and directed graphs, from bounded to linear speed, and give a sharp lower bound on the possible linear speeds in each case.
\end{abstract}

\maketitle

\section{Introduction}\label{S:intro}

A combinatorial structure $P$ consists of a (finite) set, $V(P)$ (the \textit{elements} or \textit{vertices} of $P$), and a collection of relations $\Xi(P)$ on these elements. For example, letting $V = V(G)$ and $\Xi = \{E(G)\}$, we see that any graph $G$ is a combinatorial structure in this sense. Any set of relations is permissible, but in this paper we shall only need those which define oriented edges. We say that a collection (or \textit{property}) $\P$ of combinatorial structures is \textit{hereditary} if it is closed under taking \textit{induced} sub-structures. Thus, for example, the family of graphs with no induced $C_4$ is hereditary. Write $\P_n$ for the collection of distinct (non-isomorphic) structures in a property $\P$ with $n$ vertices, and call the function $n \mapsto |\P_n|$ the \textit{speed} (or unlabelled speed) of $\P$. Also, write $\P^n$ for the collection of distinct labelled structures in $\P$ with $n$ vertices, i.e., the set of non-isomorphic pairs $(P,\phi)$, where $P \in \P_n$ and $\phi:[n] \leftrightarrow V(P)$ is a bijection (or labelling of $P$), and call the function $n \mapsto |\P^n|$ the \textit{labelled speed} of $\P$. The speed and labelled speed of a property $\P$ are both very natural measures of the ``size" of $\P$. They are also quite different from one another. For example, for the collection $\Q$ of complete bipartite graphs, we have $|\Q_n| = \lceil (n+1)/2 \rceil$ and $|\Q^n| = 2^{n-1}$, but for the collection $\R$ of paths, $|\R_n| = 1$ and $|\R^n| = n!/2$; thus the measures give different answers to the question: Which is larger, the property of being complete bipartite, or that of being a path?

We are interested in the (surprising) phenomenon that for many types of combinatorial structure, only very `few' (labelled) speeds are possible. More precisely, there often exists a family $\Fset$ of functions $f : \N \to \N$ and another function $F : \N \to \N$ with $F(n)$ {\em much} larger than $f(n)$ for every $f \in \Fset$, such that if, for each $f \in \Fset$ the speed is infinitely often larger than $f(n)$, then it is also larger than $F(n)$ for every $n \in \N$. Putting it concisely: the speed {\em jumps} from $\Fset$ to $F$.

Scheinerman and Zito~\cite{SZ} were the first to study speeds of combinatorial structures when they initiated the study of the labelled speed of hereditary graph properties. What Scheinerman and Zito showed was that the functions $n \mapsto |{\mathcal P^n}|$ are far from being `arbitrary': only certain ranges of speeds are possible. A little later, considerably stronger results were proved by Alekseev~\cite{Alekseev}, Bollob\'as and Thomason~\cite{BTbox}, \cite{BT}, \cite{BT2}, and Balogh, Bollob\'as and Weinreich~\cite{BBW1}, \cite{BBW3}, \cite{BBW4}. With hindsight, one can say that in spirit the area goes back to papers of Erd\H{o}s, Kleitman and Rothschild~\cite{EKR}, Erd\H{o}s, Frankl and R\"odl~\cite{EFR}, Kolaitis, Pr\"omel and Rothschild~\cite{KPR}, Pr\"omel and Steger~\cite{PS1}, \cite{PS2}, \cite{PS3}, \cite{PS4}, \cite{PS5}, and others; for a review of the early results, see Bollob\'as~\cite{ICM}. Later in this paper, we shall make use some of these results.

We shall mainly be interested in `low-speed' properties (those satisfying $|\P^n| \le c^n$ for some constant $c$ and sufficiently large $n$), though we shall also comment on what is known for higher speeds. In this paper we will consider properties of posets and oriented graphs; papers considering similar questions for ordered graphs and related structures include \cite{order}, \cite{KK} and \cite{Klaz}. In a forthcoming paper~\cite{tourns} we shall also consider higher speeds of hereditary properties of tournaments. Our main results are summarized in the following theorems.

Let $P = (S,<_P)$ be a partially ordered set (or poset), where $<_P$ is a partial order on the set $S$. Clearly, $P$ is a combinatorial structure. If $\P$ is a collection of posets, then let $\P_n$ and $\P^n$ be the unlabelled and labelled segments of $\P$, respectively, as described above. We begin with the unlabelled speed.

\begin{thm}\label{posets1}
If $\P$ is a hereditary property of posets, then the following assertions hold.

\begin{enumerate}
\item[$(a)$] If $|\P_n|$ is unbounded, then $|\P_n| \ge \lceil\frac{n+1}{2} \rceil$ for every $n \in \N$.\\[-1.5ex]

\item[$(b)$] If also $(|\P_n| - \lceil\frac{n+1}{2} \rceil)$ is unbounded, then $|\P_n| \ge n$ \hspace{0.02cm} $\forall\, n \in \N$.\\[-1.5ex]
\end{enumerate}
Moreover, the lower bounds are best possible.
\end{thm}

Towards the end of the paper we shall extend part $(a)$ to an arbitrary property of directed graphs (see Theorem~\ref{direct}). The proof uses the ideas of the proof of Theorem~\ref{posets1}, and also Theorem~\ref{tourn}, below.  In order to prove part $(b)$, however, we shall need a more detailed structural statement about properties of posets with speed $\left\lceil (n+1)/2 \right\rceil$, which fails to hold in the more general case.

Next, we turn to the labelled speed. Note that between two labelled structures only one isomorphism is possible, whereas in the unlabelled case there are many possible such isomorphisms. This fact makes the problem somewhat simpler, and allows us to prove the following, stronger theorem.

\begin{thm}\label{posets2}
If $\P$ is a hereditary property of posets, then one of the following assertions holds.
\begin{enumerate}
\item[$(a)$] $|\P^n| = 1$ for every $n \ge N$, for some $N \in \N$.\\[-1.5ex]

\item[$(b)$] $|\P^n|$ is a polynomial. There exists $K \in \N$ and integers $a_0,...,a_K$ (with $a_K \neq 0$), such that
$$|\P^n| = \sum_{i=0}^K a_i {n \choose i}$$ for all sufficiently large $n$. Moreover,
$$|\P^n| \ge \sum_{i=0}^K {n \choose i}$$ for all $n \ge 2K + 1$. In particular, $|\P^n| \ge n + 1$ for every $n \ge 3$.\\[-1.5ex]

\item[$(c)$] $|\P^n| \ge 2^n - 1$ for every $n \ge 6$.\\[-1.5ex]
\end{enumerate}
Moreover, each of the lower bounds is best possible.
\end{thm}

We remark that it is somewhat surprising that we are able to prove a sharp lower bound (i.e., $|\P^n| \ge \sum_{i=0}^K {n \choose i}$) for the possible speeds of order $n^K$ for every $K \in \N$; for example, we were unable in \cite{order} to prove such a sharp result for ordered graphs.

A tournament $T = (V,\psi)$ is a complete graph with an orientation on each edge. Equivalently, it is a set $S$ together with an anti-symmetric function $\psi: (S \times S) \setminus \{(x,x) : x \in S\} \to \{-1,1\}$, where we write $x \to y$ if $\psi(x,y) = 1$ and $y \to x$ otherwise.

\begin{thm}\label{tourn}
Let $\P$ be a hereditary property of tournaments. Then either
\begin{enumerate}
\item[$(a)$] $|\P_n|$ is bounded, and $\exists \,M,N \in \N$ such that $|\P_n| = M$ if $n \ge N$,\\[-1.5ex]
\end{enumerate}
\noindent or
\begin{enumerate}
\item[$(b)$] $|\P_n| \ge n-2$ for every $n \in \N$.\\[-1.5ex]
\end{enumerate}
Moreover, there exists a unique property $\P$ such that $|\P_n| =
n-2$ for every $n \in \N$ with $n \ge 4$.
\end{thm}

The structure of the paper will be as follows: in Section~\ref{poset} we discuss poset properties and prove Theorems~\ref{posets1} and \ref{posets2}, and in Section~\ref{oriented} we examine properties of various types of oriented graphs and prove Theorem~\ref{tourn}.

\section{Poset Properties}\label{poset}

Recall that a hereditary property of posets $\P$ is a collection of posets which is closed under taking induced subposets. In this section, $\P$ will be a hereditary property of posets, unless otherwise stated. Recall also that $\P_n$ is the set of unlabelled posets in $\P$ with $n$ vertices, and $\P^n$ is the set of all labelled posets in $\P$ with $n$ vertices. The reader might suggest, as another natural collection associated with $\P$, that we also consider the set of linear extensions of posets in $\P_n$ (i.e., the subset of $\P^n$ for which the labelling is monotone). However, the possible speeds in this case are just a subset of those possible for ordered graphs, since a linear extension of a poset may be thought of as an ordered graph. Thus Theorem 1 of \cite{order} determines the possible speeds up to $2^{n-1}$.

The labelled speeds of poset properties have been studied extensively, for example by Kleitman and Rothschild~\cite{KR}, who in 1975 proved that the number of all labelled posets on $n$ vertices is $2^{(1+o(1))n^2/4}$. Brightwell, Pr\"omel and Steger later gave a sharper estimate, with a simpler proof~\cite{BPS}. We shall use the Kleitman-Rothschild result to prove Theorem~\ref{poslg} below, but in the proofs of Theorems~\ref{posets1} and \ref{posets2} we shall only need the trivial upper bound $3^{{n \choose 2}}$. Other results on the speeds of certain classes of poset properties are due to (amongst others) Alon and Scheinerman \cite{ALSch}, Brightwell and Goodall \cite{BG}, and Brightwell, Grable and Pr\"omel~\cite{BGP}.

We begin by recalling the following powerful theorem about labelled speeds of hereditary properties of graphs. The results for $|\P^n| \le 2^{o(n^2)}$ are from \cite{BBW1} and \cite{BBW2}, and those for $|\P^n| \ge 2^{\Theta(n^2)}$ were proved in \cite{Alekseev}, \cite{BT} and \cite{PS5}. $B_n$ denotes the number of partitions of $[n]$.

\begin{thm}\label{BBgraphs}
Let $\P$ be a hereditary property of graphs. Then one of the following is true.
\begin{enumerate}
\item[$(a)$] $|\P^n| = \sum_{i=0}^k p_i(n) i^n$ for all sufficiently large $n$, for some $k \in \N$ and some collection
$\{p_i(n)\}_{i=0}^k$ of polynomials.\\[-1.5ex]

\item[$(b)$] $|\P^n| = n^{(1-1/k+o(1))n}$ for some $2 \le k \in \N$.\\[-1.5ex]

\item[$(c)$] $n^{(1+o(1))n} = B_n \le |\P^n| \le 2^{o(n^2)}$.\\[-1.5ex]

\item[$(d)$] $|\P^n| = 2^{(1-1/k+o(1))n^2/2}$ for some $2 \le k \in \N$.
\end{enumerate}
\end{thm}

A \textit{principal} hereditary poset family is a collection of posets not containing a fixed poset $P$ as an induced subposet. In \cite{BGP} the following theorem was proved for principal hereditary poset families. Moreover, it was shown that $|Forb(P)^n| = 2^{o(n^2)}$ if and only if $P$ has height at most two.

\begin{thm}\label{poslg}
Let $\P$ be a hereditary property of posets. Then the labelled speed of $\P$ is either $2^{(1+o(1))n^2/4}$, or $2^{o(n^2)}$.
\end{thm}

\begin{proof} Let $G_Q$ be the comparability graph of a poset $Q \in \P$. Then $\G(\P)=\{ G_Q : Q \in \P \}$ is a hereditary property of graphs. Clearly,
\begin{equation}\label{posn!}
|\G(\P)_n| \le |\P_n| \le |\G(\P)^n| \le n!\cdot |\G(\P)_n|
\end{equation}
since a poset is determined by its comparability graph and any of its linear extensions. As by Theorem \ref{BBgraphs} the (labelled) speed of $\G(\P)$ is either at least $2^{(1+o(1))n^2/4}$, or $2^{o(n^2)}$, the claim follows from (\ref{posn!}) and the Kleitman-Rothschild theorem.
\end{proof}

Theorem~\ref{poslg} concerns poset properties with close to maximal speeds; now we turn our attention to properties with low speeds. First we introduce some notation and terminology. Let $\Gamma(v) = \Gamma_G(v)$ denote the set of neighbours of a vertex $v$ in a graph $G$, and define a \textit{homogeneous block} $S$ in a graph $G$ to be a maximal set of vertices satisfying $\Gamma(x)\setminus\{y\}=\Gamma(y)\setminus\{x\}$ for all $x,y\in S$. If $G$ is the comparability graph of a poset $P$ then $P$ is said to \textit{realize} $G$. We shall write $p(v)$ for the element of $P$ corresponding to the vertex $v$ of $G$, and in the labelled cases, we shall write $p(i)$ for the poset element labelled $i$, and trust that this will not cause confusion. We shall also use $<$ to denote both $<_P$, the (partial) order in the poset, and $<_{\N}$, the usual order on the positive integers, since it will always be clear to which we refer.

We shall call the poset with no comparable elements \emph{empty}, and the poset with all elements comparable a \emph{chain}. We shall also frequently use without explanation the following trivial observations: that in a poset realizing a star, the `head' is either above all the other elements, or below all of them, and that the only poset which realizes a complete graph is a chain.

We are now ready to prove Theorems~\ref{posets1} and \ref{posets2}.

\begin{proof}[Proof of Theorem~\ref{posets1}]
As in our proof of Theorem~\ref{poslg} let $\G(\P)$ be the set of comparability graphs of posets in $\P$, and recall that $|\G(\P)_n| \le |\P_n|$.

We first show that if $|\P_n|$ is unbounded then $|\P_n| \ge \frac{n+1}{2}$ for every $n \in \N$. Consider a property of posets $\P$ such that $|\P_n|$ is unbounded, and for each graph $G \in \G(\P)$, partition $V(G)$ into homogeneous blocks $B_1, B_2, ... , B_{\ell(G)}$, let $|B_i| = t_i$, and reorder so that $t_1 \ge t_2 \ge ... \ge t_{\ell(G)}$.

Suppose first that there exist graphs $G \in \G(\P)$ with arbitrarily large values of $t_2$. If $\G(\P)$ contains all graphs consisting of a clique and an independent set and all edges in between, or the complements of all such graphs, then $|\P_n| \ge |\G(\P)_n| \ge n$, so we may assume that there exists $M \in \N$ such that whenever $t_2 \ge M$, both $B_1$ and $B_2$ are cliques, or both are independent sets. Also if $\P$ contains all posets which realize a complete bipartite graph then $|\P_n| \ge n$. But if $B_1 \cup B_2$ is a clique or an independent set, then since $B_1$ and $B_2$ are distinct homogeneous blocks, there must exist a vertex $v \in G \setminus (B_1 \cup B_2)$ which distinguishes $B_1$ from $B_2$, i.e., $B_1 \subset \Gamma(v)$ and $B_2 \subset V(G) \setminus \Gamma(v)$, or vice versa. It follows that if $B_1 \cup B_2$ is an independent set then $\G(\P)$ contains all graphs whose the edges form a star, so $|\P_n| \ge |\G(\P)_n| \ge n$, and if $B_1 \cup B_2$ induces a clique then $\G(\P)$ contains all graphs whose non-edges form a star, so again $|\P_n| \ge |\G(\P)| \ge n$. We are left only with the possibility that $B_1$ and $B_2$ are incomparable chains, and since the collection of such posets has speed $\lceil \frac{n+1}{2} \rceil$, it follows that if $t_2$ is unbounded then $|\P_n| \ge \lceil \frac{n+1}{2} \rceil$ for all $n \in \N$.

So assume that there exists $K$ such that $t_2 \le K$ for all graphs $G \in \G(\P)$. Suppose that for each $L \in \N$ there exists $v \in G \in \G(\P)$ with $L \le |\Gamma(v)| \le |G| - L$. Then by Ramsey's Theorem and the pigeonhole principle we have graphs in $\G(\P)$ with arbitrarily large values of $t_2$, contradicting our assumption. So let $L \in \N$ be such that for every $n$, each vertex $v \in G \in \G(\P)_n$ has degree at most $L$ or at least $n - L$. Given $G \in \G(\P)_n$, let $X = \{ v \in V(G) : d(v) \le L \}$, let $Y = V(G) \setminus X$. It is easy to see (by considering the number of edges between $X$ and $Y$) that $\min(|X|,|Y|) \le 2L$. Now suppose that for each $M \in \N$ there exists some $G \in \G(\P)$ with $e(G[X]) \ge M$. Then since each edge in $G[X]$ is incident with at most $2L - 2$ others, we can recursively find an arbitrarily large induced matching in some $G$. But now $\G(\P)$ contains all graphs with all degrees at most one, and since there are $\lceil \frac{n+1}{2} \rceil$ such graphs of order $n$, this implies that $|\P_n| \ge |\G(\P)| \ge \lceil \frac{n+1}{2} \rceil$. Similarly, if there are graphs $G \in \G(\P)$ with arbitrarily large values of $e(G[Y])$ then $\G(\P)$ contains all graphs with all co-degrees at most one, and we again have $|\P_n| \ge |\G(\P)| \ge \lceil \frac{n+1}{2} \rceil$.

So we may assume that there exists $M \in \N$ such that $e(G[X]) \le M$ and $e(\overline{G}[Y]) \le M$ for every $G \in \G(\P)$. It follows that for every $G \in \G(\P)_n$ we have $t_1 \ge n - (2L^2+ 2L + 2M)$.

Let $K = 2L^2 + 2(L+M)$ and consider the posets in $\P$ which realize graphs in $\G(\P)$ with $B = B_1$ an independent set. The number of such posets of order $n$ is at most $K3^{{K \choose 2}}3^K$, which is a constant.

Now consider those posets in $\P$ which realize graphs in which $B$ induces a clique (and so a chain in the poset). Suppose for each $N \in \N$ there exists a graph $G \in \G(\P)$ and a vertex $v \in V(G) \setminus B$ such that $|\{ b \in B : p(v) > p(b) \}| \ge N$ and $|\{ b \in B : p(v) < p(b) \}| \ge N$. Since $v \notin B$, there must exist another element $u \notin B$ such that either $p(u) <> p(v)$ and $p(u)$ is comparable to the elements of $B$, or $p(u) > p(v)$ say, and $p(u) <> p(b)$ for all $b \in B$. The latter case is impossible however, since $p(u) > p(v) > p(b)$ for some $b \in B$, so we must have the former, i.e., $p(u) <> p(v)$. But then for every $b \in B$ we have $p(u) > p(b)$ if and only if $p(v) > p(b)$, so $\P$ contains all the posets consisting of two incomparable elements with a chain of $N$ elements above them, and a chain of $N$ elements below them. Since $N$ was arbitrary, this in turn implies that $|\P_n| \ge n$.

We are left with the case that $t_1 \ge n - K$ in every $G \in \G(\P)_n$, $B$ induces a clique, and there exists an $N \in \N$ such that $|\{ b \in B : p(v) > p(b) \}| < N$ or $|\{ b \in B : p(v) < p(b) \}| < N$ for every $P \in \P$ and every $v \in V(G) \setminus B$. There are now only boundedly many choices for the orientations of edges from each vertex outside $B$, so the number of posets in $\P_n$ is bounded above by $K3^{{K \choose 2}}(2N+1)^K$, which contradicts our initial assumption.\\

We have proved part (a): that if $|\P_n|$ is unbounded then $|\P_n| \ge \lceil \frac{n+1}{2} \rceil$ for every $n \in \N$. Moreover, we have shown that the only possible properties with speed $\lceil \frac{n+1}{2} \rceil$ are $\Q$, the set of all posets which are the union of two (incomparable) chains, the property $\R$ consisting of posets with comparability graphs of maximum degree at most one, and the property $\overline{\R}$ consisting of all posets with comparability graphs of minimum degree $n-1$. Observe that for each such graph there exists a unique poset realizing it. It follows that in fact $|\overline{\R}_n| \ge F_n$, the $n^{th}$ Fibonacci number, since for every sequence $(a_1,...,a_s)$ with $a_i \in \{1,2\}$, $s \in \N$ and $a_1 + ... + a_s = n$, $\overline{\R}_n$ contains the (unlabelled) poset on $p_1,...,p_n$ in which $p_i < p_j$ if $i < j$ unless $a_1+...+a_t \le i < j < a_1 + ... + a_{t+1}$ for some $t \le s - 1$. Hence the only properties with speed $\lceil \frac{n+1}{2} \rceil$ are $\Q$ and $\R$, both of which have exactly this speed.

We have also shown that if $\P$ contains neither the property $\Q$, nor the property $\R$, then $|\P_n| \ge n$ for every $n \in \N$ (if $\P$ contains $\overline{\R}$ then $|\P_n| \ge F_n \ge n$). The union of $\Q$ and $\R$ has speed $ \ge n$, since $|\Q_n \cap \R_n| \le 1$ when $n \ge 3$, so if $|\P_n| < n$ for some $n$ and $|\P_n| - \lceil \frac{n+1}{2} \rceil$ is positive and unbounded then there are only a small number of possibilities.

First suppose $\P$ contains $\Q$, and for each $N \in \N$, let $\P^{(N)} = \P \setminus \{ P \in Q : $ both chains have at least $N$ elements$\}$. If $\P^{(N)}$ is hereditary for some $N \in \N$ then by the above, either its speed is bounded, or at least $n$, or $\P^{(N)}$ contains $\Q$ or $\R$. In each case we have a contradiction to one of our assumptions (if $\P^{(N)}$ contains $\R$ say, then $\Q \cup \R \subset \P$), so $\P^{(N)}$ is not hereditary for any $N \in \N$. Hence, for each $N \in \N$, there must exist $P \in \P^{(N)}$ containing an induced copy of some $Q \in \Q$ with both chains ($B_1$ and $B_2$ say) having at least $N$ elements. Take $B_1$ and $B_2$ to be maximal, subject to the condition that $P[B_1 \cup B_2] \in \Q$ (so the chains are incomparable), and observe that $B_1 \cup B_2$ does not cover $P$, since $P \in \P^{(N)}$. For such a poset $P$ in which $|B_1| = |B_2| = 2N$ say, choose a vertex $v$ of $G_P \setminus (B_1 \cup B_2)$ (where $G_P$ is the comparability graph of $P$), and consider the neighbours of $v$ in $B_1$ and $B_2$. There are various cases to consider: we shall show that in each case $|\P_N| \ge N$. Since $N$ was arbitrary, this will suffice to prove the case when $\P$ contains $\Q$.

Note first that if $|\Gamma_{B_1}(v)|, |\Gamma_{B_2}(v)| \ge N$, then $\G(\P)$ contains the graph consisting of two cliques of size $N$, with no edges in between, and one other vertex adjacent to both of the cliques. This graph has $N$ distinct subgraphs of order $N$, so $|\G(\P)_N| \ge N$. The case $|\Gamma_{B_1}(v)|, |\Gamma_{B_2}(v)| \le N$ is similar, so assume that $|\Gamma_{B_1}(v)| \ge N$ and $|\Gamma_{B_2}(v)| \le N$. Without loss of generality, assume $p(b) > p(v)$ for some $b \in B_1$.

If $|\Gamma_{B_1}(v)| = 2N$, then since $B_1$ was maximal, we must have $|\Gamma_{B_2}(v)| \ge 1$, so $b_0 \in \Gamma_{B_2}(v)$ say. Then since $p(b_0)$ is incomparable to $B_1$, we must have $p(b_0) > p(v)$, and $p(b_1) > p(v)$ for every $b_1 \in B_1$. Note also that $p(b_2) > p(v)$ for every $b_2 \in B_2$ such that $p(b_2) > p(b_0)$. Thus $\P$ contains the poset consisting of two chains of size $N$ with the top element of one above the bottom element of the other, and all other pairs of elements (not in the same chain) incomparable. This has $N$ distinct subposets of order $N$, so again $|\P_N| \ge N$.

So we may assume that $N \le |\Gamma_{B_1}(v)| < 2N$, and $|\Gamma_{B_2}(v)| \le N$, which implies that the graph consisting of a clique of size $N$ and a clique of size $N+1$ minus an edge, with no edges in between, is in $\G(\P)$. This has $N - 1 + \lceil \frac{N+1}{2} \rceil \ge N$ distinct subgraphs of order $N$, so again $|\G(\P)_N| \ge N$, and so the case $\Q \subset \P$ is complete.

The case where $\P$ contains $\R$ is similarly easy to deal with. Suppose that $\R \subset \P$, and that again $|\P_n| < n$ for some $n$ and $|\P_n| - \lceil \frac{n+1}{2} \rceil$ is positive and unbounded. Note that $\Q \not\subset \P$, as $|\R_n \cup \Q_n| \ge n$ for every $n \in \N$. From above we see that since $\Q \not\subset \P$ and $|\P_n| < n$ for some $n$, $t_2$ must be bounded. Recall that $X = \{v \in G : d(v) \le L\}$ and $Y = \{v \in G : d(v) \ge |G|-L \}$, and that (for some $L \in \N$), $V(G) = X \cup Y$ for every $G \in \G(\P)$ when $t_2$ is bounded. Since $\overline{\R} \not\subset \P$ then, again by the method above, it follows that in the comparability graphs of all but a bounded number of posets of $\P_n$ we have $|X| \ge |Y|$. Recall also that, for each $N \in \N$ and every $n$, the number of posets of order $n$ with $|X| \ge |Y|$ and $e(G[X]) \le N$ is bounded as well.

From this we may conclude that for any $N \in \N$, and for large $n$, $\P_n$ contains more than $\lceil \frac{n+1}{2} \rceil$ posets in whose comparability graphs $|X| \ge |Y|$ and $e(G[X]) \ge 2L^2N$. Since $d(v) \le L$ for $v \in X$, it follows that in such a graph there exists an induced matching on $2N$ vertices in $X$. For each $N \in \N$, let $\G^{[N]} \subset \G(\P)$ be the collection of comparability graphs in $\G(\P)$ in which $|X| \ge |Y|$ and there exists an induced matching on at least $2N$ vertices in $X$, and let $\P^{[N]}$ be the set of posets in $\P$ which realize a graph in $\G^{[N]}$. Note that $|(\P^{[N]})_n| \ge \lceil \frac{n+1}{2} \rceil + 1$ for every $N$ and every sufficiently large $n \in \N$.

Now, if $Y$ is non-empty in some $G \in \G^{[n+L]}$, then $\G(\P)$ also contains the graph consisting of $n$ triangles with one common point. This has $\lceil \frac{n-1}{2} \rceil + \lceil \frac{n+1}{2} \rceil \ge n$ distinct subgraphs of order $n$, so we may assume that for some $N_1 \in \N$, $Y$ is empty in all graphs $G \in \G^{[n]}$ with $n \ge N_1$. Also, if there exists some vertex in $G[X]$ with degree at least two in some $G \in \G^{[n + 1]}$, then $\G(\P)$ contains all graphs on $n$ vertices with one vertex of degree at least two and all others of degree at most one. There are $\lceil \frac{n-2}{2} \rceil + \lceil \frac{n+1}{2} \rceil \ge n$ such graphs, so we may assume that for some $N_2 \in \N$, $d(v) \le 1$ for every $v \in X$ in all graphs $G \in \G^{[n]}$ with $n \ge N_2$.

The proof is now complete, since there are only $\lceil \frac{n+1}{2} \rceil$ posets of order $n$ which realize a graph with $|Y| = 0$ and $d(v) \le 1$ for every $v \in X$, and this implies that $|(\P^{[N]})_n| \le \lceil \frac{n+1}{2} \rceil$ for every $N \ge \max(N_1,N_2)$. But this is a contradiction, so the case $\R \subset \P$ is complete.
\end{proof}

Next, we prove Theorem~\ref{posets2}.

\begin{proof}[Proof of Theorem~\ref{posets2}]
We shall use the following structural theorem on labelled graph properties (it is the case $k=2$ of Theorem 29 from~\cite{BBW1}). Let $\G_0$ denote the graph property in which all vertices have degree at most $1$. One can check that this is the property of smallest speed of the collection of (four) minimal properties given by the theorem. The theorem states that if $\G$ is a hereditary graph property with $|\G^n| < |\G^n_0|$ for some $n \in \N$, then $\exists$ $k,\ell \in \N$ such that all graphs $G \in \G$ can be partitioned into $t(G)+1$ parts $V(G) = A \cup B_1 \cup ... \cup B_{t(G)}$, where $t = t(G) \le \ell$ and $|A| \le k$, with each $B_i$ a homogeneous block.

Let $\P$ be a property of posets, and suppose there exists $6 \le m \in \N$ such that $|\P^m| \le 2^m - 2$. Observe that the speed of $\G_0$ satisfies $|\G^{n+1}_0| = |\G^n_0| + n|\G^{n-1}_0|$, and that $|\G^1_0| = 1$ and $|\G^2_0| = 2$, so $|\G^6_0| = 76 > 2^6 - 2$, and by induction $|\G^n_0| > 2^n - 2$ for every $n \ge 6$. It follows that in particular $|\G(\P)^m| \le |\P^m| \le 2^m - 2 < |\G^m_0|$, so we may apply the theorem to $\G(\P)$ to obtain $k$ and $\ell$.

Order the $B_i$ for each $G \in \G(\P)$ so that $|B_1| \ge ... \ge |B_t|$, and note that we may assume there exist graphs with arbitrarily large values of $|B_i|$ for each $1 \le i \le \ell$, since otherwise we may choose a larger value of $k$, a smaller value of $\ell$ and add the small $B_i$'s to $A$. Let
$$K = \max\{ k \: : \:\exists\,G \in \G(\P)\textup{ with }|G| = n\textup{ and }|A| = k\textup{ for arbitrarily large }n \},$$ and let
$$T = \max\{ t \: : \: \exists \,G \in \G(\P)\textup{ with $|G| = t$ and }|A| \ge K + 1 \}.$$

Suppose first that $\ell = 1$ and fix an $m \ge 6$ such that $|\P^m| \le 2^m - 2$. Let $B = B_1$. We first show that $|\P^n| = O(n^K)$ as $n \rightarrow \infty$. Pick a $G \in \G(\P)$ with $|G| \ge \max(m + K,T+1)$, so by the definition of $T$ we have $|A| \le K$, and hence $|B| \ge m$. Observe that $B$ is empty or complete (because it is a homogenous block), and that if it is complete then $|\P^n| \ge n!$ for all $1 \le n \le |B|$, since $\P^n$ then contains all linear orders on $[n]$. Since $m! > 2^m$ when $m \ge 4$, we may assume that $B$ is an independent set in every such $G$. But now if $v \in V(G) \setminus B$, and $p(v) > p(b)$ for some $b \in B$, then $p(v) > p(b)$ for every $b \in B$, and similarly if $p(v) < p(b)$ or $p(v) <> p(b)$ for some $b \in B$. Thus there are at most $K3^K3^{{K \choose 2}}$ unlabelled posets in $\P_n$ for every $n \ge \max(m + K,T+1)$, and each may be labelled in at most $n^K$ different ways, so $|\P^n| = O(n^K)$.

To get our lower bound on $|\P^n|$, we proceed by induction on $K$. We may again assume that $B$ is an independent set in $G$. Note that if $K = 0$ then $E(G) = \emptyset$ for every $G \in \G(\P)_n$ with $n \ge T + 1$, so $|\P^n| = 1$ for sufficiently large $n$, but not necessarily for all $n \ge 1$. If $K = 1$ however, then $\P^n$ contains the poset with $n-1$ pairwise incomparable elements, and the final element comparable to all others, which may be labelled in $n$ ways (if $n \ge 3$), and the empty poset, so $|\P^n| \ge {n \choose 1} + {n \choose 0}$ for every $n \ge 3$. So let $K \ge 2$, and observe that by the definition of $K$, there exist graphs in $\G(\P)_n$ with $|A| = K$ for all sufficiently large $n$. Each of these may be labelled in at least ${n \choose K}$ different ways if $n \ge 2K + 1$, since then $B$ is unique.

Now, if there exist graphs in $\G(\P)_n$ with $|A| = K - 1$ for all sufficiently large $n$, then we would be done by induction, so suppose not. Thus if we remove any vertex $u$ from $A$ (in any sufficiently large graph $G \in \G(\P)$), another vertex $v = v(u)$ of $A$ must fall into $B$, so in the graph $G-u$, $v$ is homogeneous to each vertex of $B$. It follows that there are no edges from $v$ to $B$ (since $B$ is independent), and either $uv \in E(G)$ and there are no edges from $u$ to $B$ in $G$, or $uv \not\in E(G)$ and all potential edges between $u$ and $B$ are in $G$ (since $v$ is not in $B$).

Partition $A$ into $A_1 = \{ u \in A : ub \in E(G)$ for every $b \in B \}$ and $A_2 = \{ u \in A : ub \not\in E(G)$ for every $b \in B \}$, and note that by the observation above, $v(u) \in A_2$ for every $u \in A$. Suppose $A_1 \neq \emptyset$ and let $u \in A_1$. Then $v(u) \in A_2$, and since $v$ falls into $B$, $vw \not\in E(G)$ for every $w \in A_2 \setminus v$. But $x = v(v(u)) \in A_2$ also, which implies that $vx \in E(G)$, a contradiction.

So $A = A_2$, which means there are no edges between $A$ and $B$ in $G$. Let $u \in A$ and take $v(u)$ as before. Observe that $uv \in E(G)$, since $u,v \in A_2$, and that $vw \not\in E(G)$ for every $w \in A \setminus \{u,v\}$ since $v$ is homogeneous to $B$ in $G-u$. Thus $\Gamma(v(u)) = \{u\}$ for each $u \in A$.

Now, applying this result to the vertex $v(u)$, we see that the only possibility for $v(v(u))$ is $u$, so also $\Gamma(u) = \{v(u)\}$. Since $u$ was arbitrary, it follows that $A$ induces a matching. Hence in all sufficiently large graphs $G \in \G(\P)$ with $|A| = K$, $E(G)$ consists of exactly $\frac{K}{2}$ independent edges. But the number of ways to partition $[K]$ into ordered pairs is $\frac{K!}{(K/2)!}$, so a poset $P \in \P_n$ realizing such a graph may be labelled in $\frac{K!}{(K/2)!}{n \choose K}$ ways. Since $\frac{K!}{(K/2)!}{n \choose K} \ge K{n \choose K} \ge \sum_{i=0}^K {n \choose i}$ if $K \ge 2$ and $n \ge 2K+1$, we are done in this case also.

Part (a) of the theorem (in the case $\ell = 1$) now follows instantly, since we have proved that if $K \ge 1$ then $|\P^n| \ge \sum_{i=0}^K {n \choose i} \ge n + 1$, whilst, as observed above, $K = 0$ implies that $E(G) = \emptyset$ for all $G \in \G(\P)_n$ with $n \ge T + 1$, so $\P^n$ contains only the poset with no comparable elements, and $|\P^n| = 1$ for $n \ge T+1$.

Note that for each $K \in \N$, the property $\Q_{(K)} = \{$bipartite posets $P = X \cup Y : |X| \le K$, $u > v$ if $u \in X, v \in Y$ and $u <> v$ if $u,v \in X$ or $u,v \in Y \}$ has labelled speed $|\Q_{(K)}^n| = \sum_{i=0}^K {n \choose i}$ for every $n \in \N$, so our lower bound is best possible.\\

We shall next prove that moreover if $\ell = 1$ then the labelled speed $|\P^n|$ is equal to a polynomial for sufficiently large $n$. The proof will once again go by induction on $K$. The case $K = 0$ is trivial from above, since $|\P^n| = 1$ for $n \ge T+1$, so let $K \ge 1$ and assume that the result is true for all smaller values of $K$. For each unlabelled poset $P \in \P_n$ with $|A| = K$ and $n \ge 2K+1$, we form the canonical poset $C = C(P)$ of $P$ as follows. First let $z = z_P : P \times P \to \{-1,0,1\}$ be the antisymmetric function such that $z(p,q) = 1$ if $p > q$, $z(p,q) = -1$ if $p < q$, and $z(p,q) = 0$ if $p = q$ or $p <> q$. Note that for any fixed $q$, $z(p,q)$ is constant for $p \in B$, because $B$ is a homogeneous set and induces an independent set in $G_P$. Now, let $C(P)$ be the poset on $K+1$ elements, with one element labelled $x$ and the other elements unlabelled, satisfying that $C(P) - p(x)$ is isomorphic to $P$ restricted to $A$ (let this isomorphism be $\psi$), and that for each vertex $v \in C(P) \setminus p(x)$, $z_{C(P)}(p(x),v) = z_P(u,\psi(v))$ for every $u \in B$. In other words, $C(P)$ is obtained from $P$ by collapsing the set $B$ into a single vertex, and labelling that vertex $x$. Since (as noted above) $z_P(u,\psi(v))$ is constant as $u$ varies over $B$, $C(P)$ always exists, and is uniquely determined.

Let $\C(\P) = \{C(P) : P \in \P_n$ for some $n \ge 2K+1$, and $|A| = K \}$, and for each $C \in \C(\P)$, let $\L(C)$ be the number of distinct ways in which one can label $C$ with $[K] \cup \{x\}$ so that $p(x)$ is labelled $x$. Note that if $n \ge 2K+1$, $\C(\P)_{n+1} \subset \C(\P)_n$, so for sufficiently large $n$, $\C(\P)_n$ is constant (and non-empty, by the definition of $K$). Thus there exist $N \in \N$ and $M \ge 1$ so that $|\C(\P)_n| = M$ when $n \ge N$. Now for some $n \ge N$, let $a_K = \sum_{C \in \C(\P)_n} \L(C)$, and observe that the number of labelled posets $P \in \P^n$ with $n \ge \max(2K+1,N)$ and $|A| = K$ is $a_K{n \choose K}$.

We now apply our induction hypothesis to the property $\widehat{\P} = \{P \in \P :$ either $|P| \le \max(2K,N-1)$, or $|A| \le K-1\}$. $\widehat{\P}$ is hereditary, has $\ell = 1$ and $\max\{ k :$ $\exists$ $G \in \G(\widehat{\P})$ with $|G| = n$ and $|A| = k$ for arbitrarily large $n \} \le K-1$, so there exist integers $a_0,...,a_{K-1}$ such that for sufficiently large $n$, $|\widehat{\P}_n| = \sum_{i=0}^{K-1} a_i{n \choose i}$. The result is now immediate.\\

So now assume that $\ell \ge 2$, again fix $m \ge 6$ satisfying $|\P^m| \le 2^m - 2$, and choose $G \in \G(\P)$ with $|B_1| \ge |B_2| \ge m$. As before, if either $B_1$ or $B_2$ is complete, then $|\P^m| \ge m!$, so assume both are empty. Also, if there is any edge between $B_1$ and $B_2$ then all edges are present (since both are homogeneous blocks). But then $|\P^m| \ge 2^m - 1$, since $\P$ contains all posets of the form $P = X \cup Y$, where $u > v$ for all $u \in X$, $v \in Y$, and all other pairs are incomparable. So we may assume that there are no edges in $G[B_1 \cup B_2]$. But since these are distinct homogeneous blocks, they must be distinguished by some vertex $v \in G \setminus (B_1 \cup B_2)$, so $B_1 \subset \Gamma(v)$ and $B_2 \subset V(G) \setminus \Gamma(v)$, say. Suppose without loss of generality that $p(v) > p(b)$ for some $b \in B_1$, so $p(v) > p(b)$ for every $b \in B_1$, since $B_1$ is an independent set. But now for every partition $X \cup Y \cup Z$ of $[m]$ with $|X| = 1$ and $|Y| \ge 1$, $\P^m$ contains the labelled poset $P(X,Y,Z)$ in which $p(i) > p(j)$ if and only if $i \in X$ and $j \in Y$. These posets are all distinct, and there are $m(2^{m-1}-1)$ of them, so $|\P^m| \ge m(2^{m-1}-1) > 2^m - 2$, a contradiction.

We conclude by noting that this lower bound is best possible, since the property $\R = \{$bipartite posets $P = X \cup Y : u > v$ if $u \in X, v \in Y$, and all other pairs incomparable$\}$ has labelled speed $|\R^n| = 2^n - 1$ for every $n \in \N$.
\end{proof}

\section{Oriented Graph Properties (OGPs)}\label{oriented}

We consider here six variants of hereditary oriented graph properties, namely unlabelled and labelled versions of oriented graphs, tournaments and directed graphs (in a directed graph edges may be both ways). A seventh variant, labelled oriented graphs in which the labelling $\phi$ is monotone, i.e., if $x \to y$ in $G$ then $\phi(x) < \phi(y)$ (this is a generalization of the monotone-labelled posets, the difference here being that transitivity is not required), will be considered in~\cite{order}, since the possible speeds in this case are a subset of those possible for ordered graphs. Notice that oriented graphs and tournaments are special cases of directed graphs, so a certain speed function is possible for OGPs and tournament properties only if it is possible for properties of directed graphs.

The range of the speeds of OGPs is between $0$ and $4^{{n\choose 2}}$, and there is a lot to explore. In this section we shall prove the existence of an initial jump in the range of realizable speeds, for each of the properties described above. More precisely, we shall show that if the speed of one of the OGPs described above is unbounded, then it is at least $\Theta(n)$. Our main task will be proving Theorem~\ref{tourn}. The other cases will either be trivial, or will reduce to the tournament case by arguments similar to those used in Section~\ref{poset}.\\

We start with a simple observation, which we shall frequently use, and which may easily proved by induction.

\begin{obs}\label{trans}
A tournament on at least $2^n$ vertices contains a transitive subtournament on at least $n$ vertices.
\end{obs}

Now, let us define four specific tournaments. Our plan in the proof of Theorem~\ref{tourn} below will be to find one of these graphs in any unbounded-speed property of tournaments.

Let $G_1^{(k)}$ be the tournament on $[2k+1]$ in which $i \rightarrow j$ if $ 1 \le i < j \le 2k$, and $i \rightarrow (2k+1)$ iff $i \ge k+1$. Consider the $k-3$ tournaments induced by taking $s$ vertices (for each $1 \le s \le k-3$) from $[k]$, $t = k-s-1$ from $[k+1,2k]$ and the vertex $2k+1$. The (out)degree sequence is $(k-2, k-3, ... , t, t, ..., s, s, ..., 2, 1)$.

Let $G_2^{(k)}$ be the tournament on $[2k+2]$ in which $i \rightarrow j$ if $1 \le i < j \le 2k+1$, and $i \rightarrow (2k+2)$ iff $i = k+1$, and consider the $k-2$ tournaments induced by taking $s$ vertices (for each $0 \le s \le k-3$) from $[k]$, $t = k-s-2$ from $[k+2,2k+1]$, and the vertices $k+1$ and $2k+2$. The degree sequence is $(k-2, k-2, k-3, ... , t+2, t+1, t+1, t-1, ... , 1, 0)$.

Let $G_3^{(k)}$ be the tournament on $[2k+2]$ in which $i \rightarrow j$ if $1 \le i < j \le 2k+1$, and $i \rightarrow (2k+2)$ iff $i \neq k+1$. Again, consider the $k-2$ tournaments induced by taking $s$ vertices (for each $0 \le s \le k-3$) from $[k]$, $t = k-s-2$ from $[k+2,2k+1]$, and the vertices $k+1$ and $2k+2$. The degree sequence is $(k-1, k-2, ... , t+2, t, t, ... , 2, 1, 1)$.

Finally, let $G_4^{(k)}$ be the tournament on $[2k+3]$ in which $i \rightarrow j$ if $1 \le i < j \le 2k+2$, and $i \rightarrow (2k+3)$ iff $i \le k$ or $i = k+2$ and consider the $k-2$ tournaments induced by taking $s$ vertices (for each $0 \le s \le k-3$) from $[k]$, $t = k-s-3$ from $[k+3,2k+2]$, and the vertices $k+1$, $k+2$ and $2k+2$. The degree sequence is $(k-1, k-2, ... , t+3, t+1, t+1, t+1, t-1, ... , 1, 0)$.

From the outdegree sequences above we can immediately distinguish all pairs of induced subtournaments, except for in $G_1^{(k)}$, where the tournaments given by $s = i$ and $s = k-i-1$ have the same sequence for each $2 \le i \le k-3$. However these can be distinguished by finding the only transitive subgraph of order $k-1$. (Note that the induced subtournaments for $s = 1$ and $s = k-2$ cannot be distinguished in this way, and in fact are isomorphic. It is for this reason that we take $s \le k-3$ in the subtournaments of $G_1^{(k)}$.) Each graph clearly also has $\overrightarrow{T_k}$ (the transitive tournament on $k$ vertices) as an induced subgraph. It follows that $G_1^{(k)}$ has at least $k-2$ distinct subgraphs of order $k$, and $G_i^{(k)}$ has at least $k-1$ distinct subgraphs of order $k$ for $i = 2,3,4$. With a little extra work it is not hard to show that for $k \ge 4$ these are in fact exactly the numbers of distinct subtournaments in each case.

Let $\P^{(i)} = \displaystyle\bigcup_k \{ H : H $ is an induced subgraph of $G_i^{(k)} \}$ for $1 \le i \le 4$.\\

Next, if $G$ is a directed graph, then given a vertex $u$ and a transitive tournament $T \subset G$, $u \not\in T$, we define the {\em pattern} $z_T(u) \in \{-1,0,1,2\}^{|T|}$ of $u$ on $T$ as follows. Let $t = |T|$, and for each $i \in [t]$, let $v_i$ be the vertex in $T$ with outdegree $t - i$ in $T$. Let $z_T(u)_i = 0$ if $u \not\to v_i$ and $v_i \not\to u$, $z_T(u)_i = 1$ if $u \to v_i$ and $v_i \not\to u$, $z_T(u)_i = -1$ if $u \not\to v_i$ but $v_i \to u$, and $z_T(u)_i = 2$ if $u \to v_i$ and $v_i \to u$. For any collection $\P$ of directed graphs, let $\cur{Z}(\P)_n$ be the set of non-transitive patterns (i.e., vectors $z$) in $\{-1,0,1,2\}^n$ which occur in $\P$, so $\cur{Z}(\P)_n = \{ z : \exists$ $u \in P \in \P$ and a transitive tournament $u \not\in T \subset P$ with $|T| = n$ such that $z_T(u) = z$, and $P[T \cup u] \neq \overrightarrow{T_{n+1}} \}$.

We first prove the following lemma, which will be needed for several of the proofs in this section.

\begin{lemma}\label{patterns}
If $\P$ is a hereditary directed graph property, and $|\cur{Z}(\P)_n|$ is unbounded, then $|\P_n| \ge n - 2$ for every $n \in \N$. Moreover, if $\P$ consists only of tournaments, then $\P^{(i)} \subset \P$ for some $i \in [4]$.
\end{lemma}

\begin{proof}
Let $k \in \N$. We shall show that $|\P_k| \ge k - 2$, and if $\P$ consists only of tournaments then $G_i^{(k)} \in \P$ for some $1 \le i \le 4$. For some large $n$, we wish to choose a collection $\G$ of pairs $(G,u)$, with each $G \in \P_{n+1} \setminus \overrightarrow{T_{n+1}}$ and $u \in V(G)$, such that for each $(G,u) \in \G$ the directed graph $T(G)$ induced by $V(G) \setminus u$ is a transitive tournament, and the patterns $z_{T(G)}(u)$ for $(G,u) \in \G$ are all distinct. Since $|\cur{Z}(\P)_n|$ is unbounded, there is an $n \in \N$ for which we can find such a $\G$ with $|\G| \ge 4^{8k}(2k+1)$.

To ease the notation, we shall write $z(G,u)$ for $z_{T(G)}(u)$, and $G_i$ for $G_i^{(k)}$ for $1 \le i \le 4$. As before, for each $i \in [n]$ let $v_i$ be the vertex in $T = T(G)$ with outdegree $n - i$ in $T$. Choose a subset $\widehat{\G} \subset \G$ of size at least $2k+1$ such that if $(G,u), (G',u') \in \widehat{\G}$, then $z(G,u)_i = z(G',u')_i$ for all $1 \le i \le 4k$ and $n-4k+1 \le i \le n$. In other words choose a subset in which $u$ has the same pattern on the `top' and `bottom' $4k$ vertices of $T(G)$. Call the remaining $n - 8k$ vertices of $T$ the `middle' vertices of $T$.

Now, let $(G,u) \in \widehat{\G}$, and for each $\ell \in \{-1,0,1,2\}$, let $A_\ell = \{ v_i \in T : 1 \le i \le 4k$ and $z(G,u)_i = \ell \}$ and $B_\ell = \{ v_i \in T : n-4k+1 \le i \le n$ and $z(G,u)_i = \ell \}$. Note that $|A_\ell|$ and $|B_\ell|$ do not depend on $(G,u)$. By the pigeonhole principle, $|A_{\ell_1}| \ge k$ and $|B_{\ell_2}| \ge k$ for some pair $\ell_1,\ell_2 \in \{-1,0,1,2\}$. For convenience later on, choose $A' \subset A_{\ell_1}$ and $B' \subset B_{\ell_2}$ with $|A'| = |B'| = k$. The remainder of the proof now consists of a fairly simple case analysis: we show that for each pair $(\ell_1,\ell_2)$, either $G_i \in \P_k$ for some $i \in [4]$, or $|\P_k| \ge k$ and $\P_k$ contains a non-tournament.

Suppose first that $\{\ell_1,\ell_2\} \not\subset \{-1,1\}$, so the directed graphs in $\widehat{\G}$ are not tournaments. If also $\ell_1 \neq \ell_2$, then for any $(G,u) \in \widehat{\G}$, the directed graph $G[u \cup A' \cup B']$ has at least $k$ distinct induced subgraphs on $k$ vertices, so $|\P_k| \ge k$. If $\ell_1 = \ell_2$ then recall that $|\widehat{\G}| \ge 2k + 1 > 1$, so for some $(G,u) \in \widehat{\G}$ and some $4k+1 \le i \le n-4k$ we have $z(G,u)_i \neq \ell_1$. Now $G[u \cup v_i \cup A' \cup B']$ has at least $k$ distinct induced subgraphs on $k$ vertices, so again $|\P_k| \ge k$.

So we may assume that $\{\ell_1,\ell_2\} \subset \{-1,1\}$, which means that $u \cup A_{\ell_1} \cup B_{\ell_2}$ induces a tournament for each $(G,u) \in \widehat{\G}$. We split into four subcases: first suppose that $\ell_1 = 1$ and $\ell_2 = -1$. This means that $u$ is `above' $k$ of the top vertices of $T$ and `below' $k$ of the bottom vertices. It is easy to see that for any $(G,u) \in \widehat{\G}$, $G[u \cup A' \cup B']$ is a copy of $G_1$. By heredity, $G_1 \in \P$.

The two cases $\ell_1 = \ell_2 \in \{-1,1\}$ are almost the same, so we shall only give the proof for $\ell_1 = \ell_2 = 1$. In this case, choose $(G,u) \in \widehat{\G}$ and $4k+1 \le j \le n-4k$ such that $z(G,u)_j \neq 1$ (this is again possible because $|\widehat{\G}| > 1$). Now if $z(G,u)_j = -1$ then $G[u \cup v_j \cup A' \cup B']$ is a copy of $G_2$, so $G_2 \in \P$, while if $z(G,u)_j \in \{0,2\}$ then $G$ is a non-tournament, and $G[u \cup v_j \cup A' \cup B']$ has $k$ distinct induced subgraphs on $k$ vertices, so $|\P_k| \ge k$. Similarly, if $\ell_1 = \ell_2 = -1$ then $G_3 \in \P$, or there exists a non-tournament in $\P$ and $|\P_k| \ge k$.

The final case, $\ell_1 = -1$ and $\ell_2 = 1$, is slightly more complicated. If there exists a $(G,u) \in \widehat{\G}$ and $4k+1 \le j \le n-4k$ such that $z(G,u)_j \notin \{-1,1\}$ (so $G$ is a non-tournament), then the directed graph $G[u \cup v_j \cup A' \cup B']$ has at least $k$ distinct induced subgraphs on $k$ vertices, and we have $|\P_k| \ge k$ as before. Thus we may assume that $z(G,u)_j \in \{-1,1\}$ for every middle vertex $v_j$ in every directed graph $G$ with $(G,u) \in \widehat{\G}$. Also, if for some $(G,u) \in \widehat{\G}$ there exists an induced cyclic triangle consisting of $u$ and two of the middle vertices of $T(G)$, then $G$ contains a copy of $G_4$ (induced by $u$, $A'$, $B'$ and these two vertices), so assume that this is not the case. It follows that for every $(G,u) \in \widehat{\G}$, $G[u \cup M]$ is a transitive tournament (where $M = \{v_i \in T(G) : 4k+1 \le i \le n-4k\}$ is the set of middle vertices).

Let $X(G,u) = \{v_i \in M : z(G,u)_i = 1\}$ and $Y(G,u) = \{v_i \in M : z(G,u)_i = -1\}$. Now, since $|\widehat{\G}| = 2k+1$, there must be some $(G,u) \in \widehat{\G}$ for which $|X(G,u)| \ge k$ and $|Y(G,u)| \ge k$. For this $G$, choose $X' \subset X(G,u)$ and $Y' \subset Y(G,u)$ with $|X'| = |Y'| = k$. Since $G \neq \overrightarrow{T_{n+1}}$, there must be some $i \in [4k]$ for which either $z(G,u)_i \neq -1$ or $z(G,u)_{n+1-i} \neq 1$. Without loss of generality suppose $z(G,u)_i \neq -1$. If $z(G,u)_i \in \{0,2\}$ then $G$ is a non-tournament and $G[u \cup v_i \cup X' \cup Y']$ has $k$ distinct subtournaments on $k$ vertices, so $|\P_k| \ge k$; if $z(G,u)_i = 1$ then $G[u \cup v_i \cup X' \cup Y'] = G_2$, so $G_2 \in \P$.

We have shown that either $|\P_k| \ge k$ (and $\P$ contains a non-tournament) or $G_i^{(k)} \in \P$ for some $i \in [4]$. Since $|\P^{(i)}_k| \ge k-2$ for each $i$, it follows that $|\P_k| \ge k-2$ in the latter case as well, and we are done.
\end{proof}

Next, we classify the bounded-speed properties of tournaments. We first define some canonical properties. If $T$ is a tournament and $A,B \subset T$, then write $A \to B$ if $a \to b$ for every $a \in A$ and $b \in B$. Given $a,b,c \in \N \cup \infty$, we say that a tournament $T$ can be $(a,b,c)$-partitioned if there exist $A$, $B$ and $C$, pairwise disjoint, with $|A| \le a$, $|B| \le b$ and $|C| \le c$, satisfying $B$ is transitive, $A \to B$ and $B \to C$. Let $\T(a,b,c)$ be the collection of tournaments $T$ which can be $(a,b,c)$-partitioned.

\begin{lemma}\label{bddt}
Let $\P$ be a hereditary property of tournaments. If $|\P_n| \le K$ for some $K \in \N$ and every $n \in \N$, then $\P \subset \T(f(K),\infty,f(K))$, where $f(K) = 2^{6K+15} + 3K + 7$.
\end{lemma}

\begin{proof}
We wish to show that for each $n$ and each tournament $G \in \P_n$, $G$ can be $(f(K),\infty,f(K))$-partitioned. Observe first that the result is trivial if $n \le 2f(K)$. Now, given a tournament $G \in \P_n$ with $n \ge 2f(K) > 2^{6K + 16}$, let $T = T(G)$ be a largest transitive subtournament in $G$, and choose a vertex $u = u(G) \in G \setminus T(G)$ (if $T(G) = G$ then $G$ is transitive and the result is trivial). Let $t = |T| \ge 6K + 16$ (by Observation~\ref{trans}), and as before, write $v_i$ for the vertex in $T$ with outdegree $t - i$ in $T$. Also let $M = \{v_i \in T : 2K+6 \le i \le t-2K-5 \}$ be the `middle' vertices of $T$ and, with foresight, define $B = \{v_i \in T : 3K+8 \le i \le t-3K-7 \}$ be the `very middle' vertices of $T$.

Consider the pattern $z = z_{T(G)}(u(G)) \in \{-1,1\}^{|T|}$ of $u$ on $T$. For $\ell \in \{-1,1\}$, let $A_\ell = \{v_i : i \in [2K+5], z_i = \ell\}$ and $B_\ell = \{v_i : t+1-i \in [2K+5], z_i = \ell \}$, as in the proof of Lemma~\ref{patterns}. Clearly, for some $\ell_1,\ell_2 \in \{-1,1\}$, $|A_{\ell_1}| \ge K+3$ and $|B_{\ell_2}| \ge K+3$. Again, choose $A' \subset A_{\ell_1}$ and $B' \subset B_{\ell_2}$ with $|A'| = |B'| = K+3$. We shall show first that in each case (i.e., for each pair $(\ell_1,\ell_2)$) either $|\P_{K+3}| \ge K+1$ (a contradiction), or $z$ is constant on $B$, the `very middle' vertices of $T$ (i.e., $z_i$ is constant for $3K+8 \le i \le t-3K-7$).

The first three cases are easy to deal with. If $\ell_1 = 1$ and $\ell_2 = -1$, then $G[u \cup A' \cup B'] = G_1^{(K+3)}$, so $G_1^{(K+3)} \in \P$, so $|\P_{K+3}| \ge K+1$. Next, if $\ell_1 = \ell_2 = 1$ then either $z_i = 1$ for every $2K+6 \le i \le t-2K-5$ (i.e. $z$ is constant on the middle vertices of $T$, and hence on the very middle vertices), or $z_j = -1$ for some $2K+6 \le j \le t-2K-5$. In the latter case, $G[u \cup v_j \cup A' \cup B'] = G_2^{(K+3)}$, so $G_2^{(K+3)} \in \P$, so $|\P_{K+3}| \ge K+1$. Similarly if $\ell_1 = \ell_2 = -1$ then either $z_i = -1$ for every $2K+6 \le i \le t-2K-5$, so $z$ is constant on the middle vertices, or $G_3^{(K+3)} \in \P$, in which case $|\P_{K+3}| \ge K+1$.

The fourth case, $\ell_1 = -1$ and $\ell_2 = 1$, is once again a little trickier. Suppose first that there exists an induced cyclic triangle consisting of $u$ and two of the middle vertices, $v_i$ and $v_j$ say, of $T(G)$. Then $G[u \cup v_i \cup v_j \cup A' \cup B'] = G_4^{(K+3)}$, so $G_4^{(K+3)} \in \P$, so $|\P_{K+3}| \ge K+1$. So we may assume that $G[u \cup M]$ is transitive.

Now, since $T$ was chosen to be maximal, $G[u \cup T]$ must be non-transitive, so there must be some $i \in [2K+5]$ for which either $z_i = 1$ or $z_{t+1-i} = -1$. Without loss of generality suppose $z_i = 1$, so $u \to v_i$.

Let $X = \{v_j \in M, z_j = 1\}$ and $Y = \{v_j \in M, z_j = -1\}$. If $|X| \ge K+3$ and $|Y| \ge K+3$ then $G[u \cup v_i \cup X \cup Y]$ contains an induced copy of $G_2^{(K+3)}$, in which case $|\P_{K+3}| \ge K+1$, so assume not. Since $G[u \cup M]$ is transitive, this implies that $z$ is constant on $B$.

We have shown that if $|\P_{K+3}| \le K$, then the vertices of $G \setminus T$ can be partitioned into sets $A' = \{v \in G \setminus T : v \to b$ for every $b \in B \}$ and $C' = \{v \in G \setminus T : b \to v$ for every $b \in B \}$. Now observe that in any $G \in \P$ we have $|A'|, |C'| \le 2^{6K+15}$, since otherwise there must exist at least $6K+15$ vertices in the larger of them which form a transitive subtournament (by Observation~\ref{trans}), and adding these vertices to $B$ we obtain a larger transitive tournament than $T$, a contradiction.

Finally let $A = A' \cup \{v_i \in T : i \le 3K+7 \} = \{u \in G : u \to B\}$, and let $C = C' \cup \{v_i \in T : i \ge t-3K-6 \} = \{v \in G : B \to v\}$. This is clearly an $(f(K),\infty,f(K))$-partition of $G$.
\end{proof}

The following lemma limits the possible oscillation of $|\P_n|$. We say that a property of tournaments $\P$ is a \textit{core} property if for every $P \in \P$, there exists $P' \in \P$, $P' \neq P$, such that $P$ is an induced subtournament of $P'$.

\begin{lemma}\label{tinc}
Let $\P$ be a core hereditary property of tournaments. If $n \ge 2m + 1$ then $|\P_n| \ge |P_m|$.
\end{lemma}

\begin{proof}
For each $G \in \P_m$, choose a $G' \in \P_N$, containing $G$ as a subtournament, with $N \ge 2^n + m$. Let the vertices corresponding to some copy of $G$ in $G'$ be $v_1,...,v_m$ (we shall abuse notation and write $G = \{v_1,...,v_m\}$), and for each $u \in G' \setminus G$ let $z(u) = \{i \in [m] : u \to v_i\}$. $z(u)$ takes only $2^m$ different values, so is constant on a set $U$ of size at least $2^{n-m}$, and so by Observation~\ref{trans}, $U$ must contain a copy $T$ of $\overrightarrow{T_{n-m}}$. For each $G \in \P_m$, let $T_n(G)$ be a tournament $G'[G \cup T] \in \P_n$ obtained in this way. We claim that the tournaments $T_n(G)$ for $G \in \P_m$ are all distinct, from which $|\P_n| \ge |P_m|$ follows immediately.

For each $G \in \P_m$, let $B = B_n(G) = G'[T]$, $A = A_n(G) = \{v \in G : v \to T\}$ and $C = C_n(G) = \{v \in G : T \to v\}$, so $T_n(G) = A \cup B \cup C$, since $z(u)$ is constant on $U$ and $T \subset U$. Notice that this is an $(m,n-m,m)$-partitioning of $T_n(G)$, and suppose that $T_n(G_1) = T_n(G_2)$ for some pair, $G_1,G_2 \in \P_m$ with $G_1 \neq G_2$. Let $f: T_n(G_1) \to T_n(G_2)$ be an edge-orientation preserving bijection, and for $i = 1,2$, let $A_i = A_n(G_i)$, $B_i = B_n(G_i)$ and $C_i = C_n(G_i)$.

We claim that $f(B_1) \subset A_2 \cup C_2$. Suppose, on the contrary, that $f(b) \in B_2$ for some $b \in B_1$. Then $f(A_1) \subset A_2 \cup B_2$ and $f(C_1) \subset B_2 \cup C_2$, since $f$ preserves edge-orientations. If $f(B_1) = B_2$, then $f$ restricted to $A_1 \cup C_1$ is an edge-preserving bijection from $G_1$ to $G_2$, a contradiction, so we may assume that $f(b') \in A_2 \cup C_2$ for some $b' \in B_1$ and $f^{-1}(b'') \in A_1 \cup C_1$ for some $b'' \in B_2$. Without loss of generality, let $f(b') \in A_2$ (otherwise reverse all the edges in all tournaments in $\P$). It follows that $f^{-1}(b'') \in C_1$, so $f(B_1) \subset A_2 \cup B_2$, again since $f$ preserves edge-orientations. Also $f(A_1) \subset A_2$ and $f^{-1}(C_2) \subset C_1$, so $|f(C_1) \cap B_2| = |f^{-1}(A_2) \cap B_1|$.

Let $s_1,...,s_k$ be the vertices of $f^{-1}(A_2) \cap B_1$, and let $t_1,...,t_k$ be the vertices of $f(C_1) \cap B_2$, ordered so that $s_i \to s_j$ and $t_i \to t_j$ if $i < j$. Define a function $g: A_1 \cup C_1 \to A_2 \cup C_2$ by $g(v) = f(v)$ if $f(v) \in A_2 \cup C_2$, and $g(v) = f(s_i)$ if $f(v) = t_i$. By the observations above, and with a little work, it can be checked that $g$ is an edge-orientation preserving bijection from $G_1$ to $G_2$, contradicting the assumption that $G_1 \neq G_2$. Hence $f(B_1) \subset A_2 \cup C_2$ as claimed.

But $f$ is a bijection, and $|B_1| = n-m \ge m+1 > |G_2| = |A_2 \cup C_2|$, so $f(B_1)$ cannot be a subset of $A_2 \cup C_2$. It follows that our initial assumption was false, and the tournaments $T_n(G)$ with $G \in \P_m$ are all distinct. Since $T_n(G) \in \P_n$ for every $G \in \P_m$, we have $|\P_n| \ge |P_m|$.
\end{proof}

We are now ready to prove the main result of the section.

\begin{lemma}\label{tourns}
If $\P$ is a hereditary property of tournaments, and $|\P_n|$ is unbounded, then $\P^{(i)} \subset \P$ for some $1 \le i \le 4$.
\end{lemma}

\begin{proof}
Let $k \in \N$ and $\P$ be a hereditary property of tournaments with $|\P_n|$ unbounded, and $G_i^{(k)} \notin \P$ for each $i \in [4]$. We shall write $G_i$ for $G_i^{(k)}$.

Suppose first that there exists a constant $K$ such that from each tournament $G \in \P$ we may remove at most $K$ vertices and leave a transitive tournament. For each $G \in \P$, choose a maximal transitive tournament $T = T(G)$ on at least $|G|-K$ vertices, let $t = |T|$, and again let $v_i$ be the vertex of $T$ with outdegree $t - i$ in $T$.

Recall that $\cur{Z}(\P)_n = \{ z \in \{-1,1\}^n : \exists$ $u \in G \in \P$ and a transitive tournament $u \not\in T \subset G$ with $|T| = n$, such that $z_T(u) = z$, and $G[T \cup u] \neq \overrightarrow{T_{n+1}} \}$. If $|\cur{Z}(\P)_n| \le L$ for some $L \in \N$ and every $n \in \N$, then $|\P_n| \le (K+1)2^{K^2}L^K$, a contradiction. (This inequality follows because there are fewer than $(K+1)2^{K^2}$ ways to choose the tournament induced by $G \setminus T$, and at most $L^K$ ways to orient the edges between $T$ and $G \setminus T$ if $G \neq \overrightarrow{T_n}$.) But if $|\cur{Z}(\P)_n|$ is unbounded then, by Lemma~\ref{patterns}, $\P^{(i)} \subset \P$ for some $i \in [4]$, which also contradicts our initial assumptions.

So we may assume that for each integer $K$, there exists $G \in \P$ containing no transitive subtournament on $|G|-K$ vertices. Let $K = 4k {{2k} \choose k} 2^{k+1}$, and choose such a $G$ and a maximal transitive subtournament $T = T(G)$ of $G$.

Call two vertices of $T$ `adjacent' if they have the same orientation with respect to the other vertices in $T$, and note that by maximality of $T$, each vertex of $G \setminus T$ must form a cyclic triangle with some two adjacent vertices of $T$. For each vertex $u \in G \setminus T$, choose such a cyclic triangle, $(u,v_j,v_{j+1})$ say, and observe that $v_j$ and $v_{j+1}$ must lie in the top or in the bottom $2k+1$ vertices of $T$ (i.e., $j \in [2k] \cup [t-2k,t-1]$), else for some $1 \le i \le 4$ there would be a copy of $G_i$ in $\P$. (This follows by the now-familiar method: choose $X \subset \{v_1,...,v_{2k-1}\}$ and $Y \subset \{v_{t-2k+2},...,v_t\}$ with $|X| = |Y| = k$ and $z_T(u)_i = \ell_1$ for every $v_i \in X$, $z_T(u)_j = \ell_2$ for every $v_j \in Y$, and consider the four cases $(\ell_1,\ell_2)$ with $\ell_1,\ell_2 \in \{-1,1\}$.)

Now, by our choice of $K$, there must exist at least ${{2k} \choose k} 2^{k+1}$ vertices $u \in G \setminus T$ for which we chose the same pair of vertices $\{v_j,v_{j+1}\} \subset T$. Let $U = \{ u \in G \setminus T : u \to v_j$ and $v_{j+1} \to u \}$. Suppose, without loss of generality, that this pair lies in the top $2k+1$ vertices (so $j \in [2k]$), and for each $u \in U$ choose $k$ vertices $w_1,...,w_k$ from the bottom $2k$ of $T$ which have the same orientation with respect to $u$. For at least $2^{k+1}$ of the vertices in $U$ we chose the same $k$ vertices (call these $k$ vertices $B$), and for at least $2^k$ of these, the orientation of $uw_i$ is the same. By Observation~\ref{trans}, these $2^k$ vertices must contain a transitive tournament on at least $k$ vertices.

Let $A$ be this set of $k$ vertices from $U$, let $B$ be the $k$-set defined above, and let $C = \{v_j,v_{j+1}\}$. Rename $v_j$ and $v_{j+1}$ $c_1$ and $c_2$ respectively, and recall that for two sets of vertices, $X$ and $Y$, we write $X \to Y$ if $x \to y$ for every $x \in X$ and $y \in Y$. Thus we have a tournament in $\P$ on $2k+2$ vertices consisting of three disjoint transitive subtournaments $A$, $B$ and $C$, with $|A| = |B| = k$ and $|C| = 2$, and with either $A \to B$ or $B\to A$; $B \to C$ or $C \to B$; and with $A \to c_1$, $c_1 \to c_2$ and $c_2 \to A$.

Now if $A \rightarrow B$ and $B \rightarrow C$, remove $c_1$, and if $B \rightarrow A$ and $C \rightarrow B$ remove $c_2$; in each case we get a copy of $G_1$. If $A \rightarrow B$ and $C \rightarrow B$ the tournament is $G_2$; if $B \rightarrow A$ and $B \rightarrow C$ the tournament is $G_3$. In each case we have shown that $G_i \in \P$ for some $1 \le i \le 4$, and we have the desired contradiction.
\end{proof}

The proof of Theorem~\ref{tourn} now follows almost immediately from Lemmas~\ref{bddt}, \ref{tinc} and \ref{tourns}.

\begin{proof}[Proof of Theorem~\ref{tourn}]
If $|\P_n| \le n-3$ for some $n \in \N$, then $\P^{(i)} \not\subset \P$ for any $i \in [4]$, so by Lemma~\ref{tourns}, $|\P_n|$ is bounded. Now, by Lemma~\ref{bddt} we can $(f(K),\infty,f(K))$-partition every tournament $T \in \P$, where $K = \max\{|\P_n| : n \in \N\}$. If $n \ge 2^{6K+18} > 4f(K) + 1$, then in such a partition the middle (transitive) block $B$ has order $|B| > (n + 1)/2$. Remove a vertex from the middle block of each $T \in \P_n$, to get a tournament $T' \in \P_{n-1}$. Suppose $T'_1 = T'_2$ for two tournaments formed in this way from $T_1,T_2 \in \P_n$, and that $f': T'_1 \to T'_2$ is an edge-orientation preserving bijection. Since, $|B'_1| = |B'_2| = |B| - 1 > (n-1)/2$, $f'(b) \in B'_2$ for some $b \in B'_1$. Now, by the proof of Lemma~\ref{tinc}, it can be seen that $f'$ extends in the obvious way to an edge-orientation preserving bijection $f : T_1 \to T_2$, so $T_1 = T_2$

It follows that $|\P_n|$ is decreasing if $n \ge 2^{6K+18}$ and so $|\P_n|$ is constant for sufficiently large $n$. Moreover by Lemma~\ref{tinc}, in a core property this constant is the maximal value taken by the speed, i.e., $|\P_n| = K$ for sufficiently large $n$.

Finally, if $|\P_n|$ is unbounded, then by Lemma~\ref{tourns}, $\P^{(i)} \subset \P$ for some $i$. Since $\P^{(i)}_n = n-1$ for large $n$ if $i > 1$, it follows that $\P^{(1)}$ is the unique property $\P$ with $|\P_n| = n -2$ for every $n \in \N$ with $n \ge 4$.
\end{proof}

\hspace{0cm}\\We finish by shortly sketching proofs for the other types of OGP.\\

Labelled tournament properties: By Observation~\ref{trans}, any tournament on $N \ge 2^n$ vertices contains a transitive tournament on $n$ vertices. Since such a transitive tournament may be labelled in $n!$ ways, this gives that the speed is either $0$ or at least $n!$.\\

Unlabelled directed graph properties: The following theorem generalizes Theorem~\ref{posets1}$(a)$ to an arbitrary hereditary property of directed graphs, but does not attempt to lay the groundwork necessary to prove a statement corresponding to part $(b)$. The proof uses several of the results above.

\begin{thm}\label{direct}
If $\P$ is a hereditary property of directed graphs and $|\P_n|$ is unbounded then $|\P_n| \ge \lceil \frac{n+1}{2} \rceil$ for every $n \in \N$.
\end{thm}

\begin{proof}
We shall consider three different graph properties induced by $\P$, and use the fact that if the speed of any of them is unbounded, then it is at least $\lceil \frac{n+1}{2} \rceil$.

Given a directed graph $D$, let the single-edge graph $K_D$ of $D$ be the graph on $V(D)$ with edges corresponding to single edges of $D$. In other words, $uv \in E(K_D)$ if $u \to v$ but $v \not\to u$, or $v \to u$ and $u \not\to v$, and $uv \notin E(K_D)$ otherwise. Let the double-edge graph $L_D$ of $D$ be the graph with edges corresponding to double edges of $D$, and $M_D$ be the graph with edges corresponding to non-edges of $D$. Let $\K = \{ K_D : D \in \P \}$, $\L = \{ L_D : D \in \P \}$ and $\M = \{ M_D : D \in \P \}$. $\K$, $\L$ and $\M$ are clearly hereditary, and $|\P_n| \ge |\K_n|, |\L_n|, |\M_n|$ for every $n \in \N$.

We use the following theorem about graph properties, which is from \cite{BBSS}, but can also be read out of the proof of Theorem~\ref{posets1}: the unlabelled speed of a hereditary property of graphs $\G$ satisfies either $|\G_n| \ge \lceil \frac{n+1}{2} \rceil$ for every $n \in \N$, or $|\G_n|$ is bounded, in which case there exists a constant $k \in \N$ such that every graph $G \in \G_n$ contains a homogeneous set $H(G)$ of size at least $n - k$. Applying this to $\K$, $\L$ and $\M$, we see that we may assume that for some $k \in \N$ (the largest of the three values given by the theorem), each graph $G \in \K_n \cup \L_n \cup \M_n$ contains a homogeneous set $H(G)$ of size at least $n - k$.

Now, let $D$ be any directed graph in $\P_n$, and define $H(D)$ to be the intersection $H(K_D) \cap H(L_D) \cap H(M_D)$. This has size at least $n - 3k$, and is a homogeneous set in single edges, double edges and non-edges. This implies that $H(D)$ either induces a complete double-edge graph, a tournament, or an independent set in $D$.

Consider first $\P^{(1)} \subset \P$ consisting of those directed graphs $D$ in $\P$ for which $H(D)$ is an independent set or a double-edge clique. Since $H(D)$ is a homogeneous set in $K_D$, $L_D$ and $M_D$, each vertex $u \in D \setminus H(D)$ sends either all double edges, all single edges, or all non-edges to $H(D)$. We shall show that the number of such graphs of order $n$ is either bounded or at least $n+1$.

For $D \in \P$ and $u \in D \setminus H(D)$, define $X(u) = \{v \in H(D) : u \to v$ and $v \not\to u \}$ and $Y(u) = \{v \in H(D) : v \to u$ and $u \not\to v \}$. Suppose that for every $n \in \N$, there exists some vertex $u$ and some graph $D$ with $u \in D \setminus H(D)$, such that both $|X(u)| \ge n$ and $|Y(u)| \ge n$ (so $u$ sends all single edges into $H(D)$). Then for each $n$, the directed graph $D[u \cup X(u) \cup Y(u)]$ has $n+1$ distinct subgraphs on $n$ vertices, so $|\P^{(1)}_n| \ge n+1$ for every $n \in \N$.

Thus we may assume that there is an $N \in \N$ such that for each such $D$, and each $u \in D \setminus H(D)$, either $|X(u)| < N$ or $|Y| < N$. But now $|\P^{(1)}_n|$ is bounded, since there are only at most $(3k+1)4^{(3k)^2}$ ways to choose $D[V(D) \setminus H(D)]$, at most $(2N+2)^{3k}$ ways to choose the cross-edges, and two ways to choose the edges inside $H(D)$.

So suppose there are an unbounded number of directed graphs $D \in \P_n$ for which $H(D)$ induces a tournament in $D$, and (for each $n \in \N$) let $\P^{(2)}_n = \P_n \setminus \P^{(1)}_n$ be the collection of such graphs. Since every tournament induced by some $H(D)$ is in $\P$, by Theorem~\ref{tourn} we have that either $|\P^{(2)}_n| \ge n-2$ for every $n \in \N$, or there are only a bounded number of different tournaments obtained in this way.

In the former case we are done, so let us assume the latter holds. Now apply Lemma~\ref{bddt}. This tells us that for some $\ell \in \N$, every tournament induced by some $H(D)$ may be $(\ell,\infty,\ell)$-partitioned. The central part of this partition is a transitive tournament on at least $|D| - 3k - 2\ell$ vertices, so an unbounded number of directed graphs $D \in \P_n$ contain such a transitive tournament. For each $D \in \P$, let $T(D)$ be the largest transitive tournament contained in $D$.

Now, since $\{|D \setminus T(D)| : D \in \P^{(2)} \}$ is bounded, and all patterns on $T(D)$ are non-transitive, there must be an unbounded number of patterns (on transitive tournaments of order $n$) which occur in $\P^{(2)}$, and hence in $\P$. The result now follows by Lemma~\ref{patterns}.
\end{proof}

Observe that this lower bound is sharp, since it is achieved by the unlabelled speed of a hereditary graph property. Two extremal examples are the property consisting of all directed graphs which are union of two double-edge cliques, and the property consisting of all directed graphs which are the union of two transitive tournaments.\\

Labelled directed graph properties: If $\P$ is a hereditary property of directed graphs, and $|\P^n| \ge 3$ for infinitely many $n$, then $|\P^n| \ge n + 1$ for every $n \ge 3$.

\begin{proof}
Form the comparability graph properties $\K$, $\L$ and $\M$ as in the unlabelled case, and observe that $|\P^n| \ge |\K^n|, |\L^n|, |\M^n|$ for every $n$. We use the following simple theorem about graph properties from~\cite{BBW1}: for any graph property $\G$, either $|\G^n| \ge n+1$ for every $n \ge 3$, or $\G_n \subset \{E_n,K_n\}$ for every $n \in \N$.

Applying this result to each of $\K$, $\L$ and $\M$, we see that if $|\P^n| \le n$ for some $n \ge 3$, then $\P_n \subset \{E_n, DK_n\} \cup \{$tournaments on $n$ vertices$\}$, where $DK_n$ denotes the complete double-edged graph on $n$ vertices. Since (by Observation~\ref{trans}) any tournament on $N \ge 2^n$ vertices contains a transitive tournament on $n$ vertices, which may be labelled in $n!$ different ways, it follows that for sufficiently large $n$, $\P_n \subset \{E_n, DK_n\}$, so $|\P^n| \le 2$.
\end{proof}

Observe again that the lower bound is sharp, since the property consisting of all double-edged stars (the empty graph is included), has speed exactly $n+1$ for $n \ge 3$, as does the property consisting of all single-edged stars, with the edges all directed away (say) from the centre.\\

Unlabelled and labelled oriented graph properties: Since an oriented graph is a special type of directed graph, the possible speeds are a subset of those possible for directed graphs, so we have the same results as for directed graphs. By the examples given above, the lower bounds are once again best possible.

\end{document}